\documentclass[12pt]{article}
\usepackage{amsthm,color,amsmath,graphicx,mathptmx,amssymb}
\usepackage{graphicx,anysize}
\usepackage[colorlinks,urlcolor=blue,linkcolor=blue,citecolor=blue]{hyperref}
\marginsize{25mm}{25mm}{15mm}{15mm}

\newtheorem{thm}{Theorem}[section]
\newtheorem{prop}[thm]{Proposition}
\newtheorem{lem}[thm]{Lemma}
{\bf}{\it}

\theoremstyle{remark}
\newtheorem{rem}[thm]{Remark}
\newtheorem*{pf}{Proof}
\theoremstyle{definition}
\newtheorem{defn}[thm]{Definition}

\allowdisplaybreaks
\numberwithin{equation}{section}

\begin{document}

\title{Standing waves with large frequency for $4$-superlinear Schr\"{o}dinger-Poisson systems
\thanks{Supported by NSFC (11171204).}
}

\author{Huayang Chen $^{\text{a}}$\and  Shibo Liu $^{\text{b}}$\thanks{email: liusb@xmu.edu.cn}
}

\date{\small\it $^{\text{a}}$ Department of Mathematics, Shantou University, Shantou 515063, China\\
$^{\text{b}}$ School of Mathematical Sciences, Xiamen University, Xiamen 361005, China}
\maketitle

\begin{abstract}
\noindent We consider standing waves with frequency $\omega$ for $4$-superlinear Schr\"{o}dinger-Poisson system. For large $\omega$ the problem reduces to a system of elliptic equations in $\mathsf{R}^3$ with potential indefinite in sign. The variational functional does not satisfy the mountain pass geometry. The nonlinearity considered here satisfies a condition which is much weaker than the classical (AR) condition and the condition (Je) of Jeanjean. We obtain nontrivial solution and, in case of odd nonlinearity an unbounded sequence of solutions via the local linking theorem and the fountain theorem, respectively.\medskip

\smallskip

\noindent\emph{MSC20000: }58E05; 35J60
\end{abstract}

\section{Introduction}

When we are looking for standing wave solutions $\psi(t,x)=\mathrm{e}%
^{-i\omega t/\hbar}u(x)$ for the nonlinear Schr\"{o}dinger equation%
\begin{equation}
i\hbar\frac{\partial\psi}{\partial t}=-\frac{\hbar^{2}}{2m}\Delta\psi
+U(x)\psi+\phi\psi-\tilde{g}(\left\vert \psi\right\vert )\psi\text{,\qquad
}(t,x)\in\mathsf{R}^{+}\times\mathsf{R}^{3} \label{e0}%
\end{equation}
coupled with the Poisson equation%
\[
-\Delta\phi=\left\vert \psi\right\vert ^{2}\text{,\qquad}x\in\mathsf{R}%
^{3}\text{,}%
\]
we are led to a system of elliptic equations in $\mathsf{R}^{3}$ of the form%
\begin{equation}
\left\{
\begin{array}
[c]{ll}%
-\Delta u+V(x)u+\phi u=g(u)\text{,} & \text{in }\mathsf{R}^{3}\text{,}\\
-\Delta\phi=u^{2}\text{,} & \text{in }\mathsf{R}^{3}\text{,}%
\end{array}
\right.  \label{e1}%
\end{equation}
where the potential $V(x)=U(x)-\omega$ and, without lose of generality we
assume $\hbar^{2}=2m$, so that the coefficient of $\Delta u$ in the first
equation is $-1$. Due to the physical context, the nonlinearity $g(t)=\tilde
{g}(\left\vert t\right\vert )t$ satisfies%
\begin{equation}
\lim_{\left\vert t\right\vert \rightarrow0}\frac{g(t)}{t}=0\text{.} \label{e}%
\end{equation}

The problem \eqref{e1} has a variational structure. It is known that there is
an energy functional $\mathcal{J}$ on $H^{1}( \mathsf{R}^{3}) \times
\mathcal{D}^{1,2}( \mathsf{R}^{3})$,
\[
\mathcal{J}(u,\phi)=\frac{1}{2}\int_{\mathsf{R}^{3}}\left(  \left\vert \nabla
u\right\vert ^{2}+V(x)u^{2}\right)  \mathrm{d}x-\frac{1}{4}\int_{\mathsf{R}%
^{3}}\left\vert \nabla\phi\right\vert ^{2}\mathrm{d}x+\frac{1}{2}%
\int_{\mathsf{R}^{3}}\phi u^{2}\mathrm{d}x-\int_{\mathsf{R}^{3}}%
G(u)\mathrm{d}x\text{,}
\]
such that $( u,\phi) $ solves \eqref{e1} if and only if it is a critical point
of $\mathcal{J}$. However, the functional $\mathcal{J}$ is strongly indefinite
and is difficult to investigate.

For $u\in H^{1}(\mathsf{R}^{3})$, it is well known that the Poisson equation%
\[
-\Delta\phi=u^{2}%
\]
has a unique solution $\phi=\phi_{u}$ in $\mathcal{D}^{1,2}(\mathsf{R}^{3})$.
Let%
\begin{equation}
\Phi(u)=\frac{1}{2}\int_{\mathsf{R}^{3}}\left(  \left\vert \nabla u\right\vert
^{2}+V(x)u^{2}\right)  \mathrm{d}x+\frac{1}{4}\int_{\mathsf{R}^{3}}\phi
_{u}u^{2}\mathrm{d}x-\int_{\mathsf{R}^{3}}G(u)\mathrm{d}x\text{,} \label{e3}%
\end{equation}
where%
\[
G(t)=\int_{0}^{t}g(\tau)\mathrm{d}\tau\text{.}%
\]
It is well known that under suitable assumptions, $\Phi$ is of class $C^{1}$
on some Sobolev space and, if $u$ is a critical point of $\Phi$, then
$(u,\phi_{u})$ is a solution of \eqref{e1}; see e.g. \cite[pp. 4929--4931]%
{MR2548724} for the details. In other words, finding critical points
$(u,\phi)$ of $\mathcal{J}$ has been reduced to looking for critical points
$u$ of $\Phi$. The idea of this reduction method is originally due to Benci
and his collaborators \cite{MR1659454,MR1714281}.

In this paper we assume that the potential $V$ and the nonlinearity $g$
satisfy the following conditions.

\begin{enumerate}
\item[$( V) $] $V\in C( \mathsf{R}^{3}) $ is bounded from below and, $\mu(
V^{-1}(-\infty,M]) <\infty$ for every $M>0$, where $\mu$ is the Lebesgue
measure on $\mathsf{R}^{3}$.

\item[$( g_{0}) $] $g\in C( \mathsf{R}) $ satisfies \eqref{e} and there exist
$C>0$ and $p\in( 4,6) $ such that%
\[
\left\vert g( t) \right\vert \leq C( 1+\left\vert t\right\vert ^{p-1})
\text{.}%
\]

\item[$(g_{1})$] there exists $b>0$ such that $4G(t)$ $\leq tg(t)+bt^{2}$, and%
\begin{equation}
\lim_{\left\vert t\right\vert \rightarrow\infty}\frac{g(t)}{t^{3}}%
=+\infty\text{.} \label{e2}%
\end{equation}

\end{enumerate}

We emphasize that unlike all previous results about the system \eqref{e1}, see
e.g. \cite{MR2769159,MR2422637,MR2548724,MR2318269,MR2428280}, we have not
assumed that the potential $V$ is positive. This means that we are looking for
standing waves of \eqref{e0} with \emph{large} frequency $\omega$.

When $V$ is positive, the quadratic part of the functional $\Phi$ is
positively definite, and $\Phi$ has a mountain pass geometry. Therefore, the
mountain pass lemma \cite{MR0370183} can be applied. In our case, the
quadratic part may posses a nontrivial negative space $X^{-}$, so $\Phi$ no
longer possesses the mountain pass geometry. A natural idea is that we may try
to apply the linking theorem (also called generalized mountain pass theorem)
\cite[Theorem 2.12]{MR1400007}. Unfortunately, due to the presence of the term
involving $\phi_{u}$, it turns out that $\Phi$ does not satisfy the required
linking geometry either. To overcome this difficulty, we will employ the idea
of local linking \cite{MR1312028,MR2135818}.

Since the term involving $\phi_{u}$ in the expression of $\Phi$ is homogeneous
of degree $4$ (see \eqref{phi}), it is natural to consider the case that
\eqref{e2} holds. In this case, we say that the nonlinearity $g( t) $ is
$4$-superlinear. In the study of such problems, the following
Ambrosetti--Rabinowitz condition \cite{MR0370183}

\begin{enumerate}
\item[(AR)] there exists $\theta>4$ such that $0<\theta G( t) \leq tg( t) $
for $t\neq0$
\end{enumerate}
is widely used, see \cite{MR2548724,MR2428280}. Another widely used condition
is the following condition introduced by Jeanjean \cite{MR1718530}

\begin{enumerate}
\item[(Je)] there exists $\theta\geq1$ such that $\theta\mathcal{G}( t)
\geq\mathcal{G}( st) $ for all $s\in\left[  0,1\right]  $ and $t\in\mathsf{R}%
$, where $\mathcal{G}( t) =g( t) t-4G( t) $.
\end{enumerate}
It is well known that if $t\mapsto\left\vert t\right\vert ^{-3}g( t) $ is
nondecreasing in $( -\infty,0) $ and $( 0,\infty) $, then (Je) holds.
Obviously, our condition $( g_{1}) $ is weaker than both (AR) and (Je).
Therefore, it is interesting to consider $4$-superlinear problems under the
condition $( g_{1}) $.

The condition $( g_{1}) $ is motivated by Alves \emph{et.\ al}
\cite{MR2769159}. Assuming in addition%
\begin{equation}
\alpha=\inf_{\mathsf{R}^{3}}V>0 \label{v}%
\end{equation}
and $b\in\lbrack0,\alpha)$, they were able to show that all Cerami sequences
of $\Phi$ are bounded. In the present paper we will show that with the compact
embedding $X\hookrightarrow L^{2}( \mathsf{R}^{3}) $ mentioned below, we can
get the boundedness of Palais-Smale sequences under much weaker condition $(
g_{1}) $, see Lemma \ref{lem1}.

Since $V$ is bounded from below, we may chose $m>0$ such that%
\[
\tilde{V}(x):=V(x)+m>1\text{,\qquad for all }x\in\mathsf{R}^{3}\text{.}%
\]
A main difficulty to solve problem \eqref{e1} is that the Sobolev embedding
$H^{1}(\mathsf{R}^{3})\hookrightarrow L^{2}(\mathsf{R}^{3})$ is not compact.
Thanks to the condition $(V)$, this difficulty can be overcame by the compact
embedding $X\hookrightarrow L^{2}(\mathsf{R}^{3})$ of Bartsch and Wang
\cite{MR1349229}, where%
\[
X=\left\{  u\in H^{1}(\mathsf{R}^{3})\left\vert \,\int_{\mathsf{R}^{3}%
}V(x)u^{2}\mathrm{d}x<\infty\right.  \right\}
\]
is a linear subspace of $H^{1}(\mathsf{R}^{3})$, equipped with the inner
product%
\[
\left\langle u,v\right\rangle =\int_{\mathsf{R}^{3}}\left(  \nabla
u\cdot\nabla v+\tilde{V}(x)uv\right)  \mathrm{d}x
\]
and the corresponding norm $\left\Vert u\right\Vert =\left\langle
u,u\right\rangle ^{1/2}$. With this inner product, $X$ is a Hilbert space. We
also note that if $V$ is coercive, namely%
\[
\lim_{\left\vert x\right\vert \rightarrow\infty}V(x)=+\infty\text{,}%
\]
then $(V)$ is satisfied.

Under our assumptions, the functional $\Phi$ given in \eqref{e3} is of class
$C^{1}$ on $X$ and, to solve \eqref{e1} it suffices to find critical points of
$\Phi\in C^{1}(X)$.

According to the compact embedding $X\hookrightarrow L^{2}( \mathsf{R}^{3}) $
and the spectral theory of self-adjoint compact operators, it is easy to see
that the eigenvalue problem%
\begin{equation}
-\Delta u+V( x) u=\lambda u\text{,\qquad}u\in X \label{eig}%
\end{equation}
possesses a complete sequence of eigenvalues%
\[
-\infty<\lambda_{1}\leq\lambda_{2}\leq\lambda_{3}\leq\cdots\text{,\qquad
}\lambda_{k}\rightarrow+\infty\text{.}%
\]
Each $\lambda_{k}$ has been repeated in the sequence according to its finite
multiplicity. We denote by $\phi_{k}$ the eigenfunction of $\lambda_{k}$, with
$\left\vert \phi_{k}\right\vert _{2}=1$, where $\vert\cdot\vert_{q}$ is the
$L^{q}$ norm. Note that the negative space $X^{-}$ of the quadratic part of
$\Phi$ is nontrivial, if and only if some $\lambda_{k}$ is negative. Actually,
$X^{-}$ is spanned by the eigenfunctions corresponding to negative
eigenvalues.\medskip

We are now ready to state our results.

\begin{thm}
\label{t1}Suppose $( V) $, $( g_{0}) $ and $( g_{1}) $ are satisfied. If $0$
is not an eigenvalue of \eqref{eig}, then the Schr\"{o}dinger-Poisson system
\eqref{e1} has at least one nontrivial solution $( u,\phi) \in X\times
\mathcal{D}^{1,2}( \mathsf{R}^{3}) $.
\end{thm}

As far as we know, this is the first existence result for \eqref{e1} in the
case that the Schr\"{o}dinger operator $S=-\Delta+V$ is not necessary
positively definite. In case that $g$ is odd, we can obtain an unbounded
sequences of solutions.

\begin{thm}
\label{t2}If $( V) $, $( g_{0}) $, $( g_{1}) $ are satisfied, and $g$ is odd,
then the Schr\"{o}dinger-Poisson system \eqref{e1} has a sequence of solutions
$( u_{n},\phi_{n}) \in X\times\mathcal{D}^{1,2}( \mathsf{R}^{3}) $ such that
the energy $\mathcal{J}(u_{n},\phi_{n})\to+\infty$.
\end{thm}

\begin{rem}
\label{re13}

\begin{enumerate}
\item[(i)] Note that if $u$ is a critical point of $\Phi$ then $\Phi
(u)=\mathcal{J}(u,\phi_{u})$. Therefore, to prove Theorem \ref{t2} it suffices
to find a sequence of critical points $\{u_{n}\}$ of $\Phi$ such that
$\Phi(u_{n})\to+\infty$.

\item[(ii)] Theorem \ref{t2} improves the recent results in \cite{MR2548724}
and \cite{MR2606783}. In these two papers the authors assumed in addition
\eqref{v}, and (AR) or (Je) respectively.
\end{enumerate}
\end{rem}

Schr\"{o}dinger-Poisson systems of the form \eqref{e1} has been extensively
studied in recent years. To overcome the difficulty that the embedding
$H^{1}(\mathsf{R}^{3})\hookrightarrow L^{2}(\mathsf{R}^{3})$ is not compact,
many authors restrict their study to the case that $V$ is radially symmetric,
or even a positive constant, see e.g. \cite{MR2417922,MR2099569,MR2230354}.
Obviously, in this case replacing $X$ with the radial Sobolev space
$H^{1}_{\mathrm{rad}(\mathsf{R}^{3})}$, the conclusions of Theorems \ref{t1}
and \ref{t2} remain valid.

\section{Tools from critical point theory}

Evidently, the properties of $\phi_{u}$ play an important role in the study of
$\Phi$. According to \cite[Theorem 2.2.1]{MR1625845} we know that%
\begin{equation}
\phi_{u}(x)=\frac{1}{4\pi}\int_{\mathsf{R}^{3}}\frac{u^{2}(y)}{\left\vert
x-y\right\vert }\mathrm{d}y\text{.} \label{phi}%
\end{equation}
Using this expression, a more complete list of properties can be found in
\cite[Lemma 1.1]{MR2769159}. Here, we only recall the ones that will be used
in our argument.

\begin{prop}
\label{p1}There is a positive constant $a_{1}>0$ such that for all $u\in X$ we
have%
\begin{equation}
0\leq\frac14\int_{\mathsf{R}^{3}}\phi_{u}u^{2}\mathrm{d}x\leq a_{1}\left\Vert
u\right\Vert ^{4}\text{.} \label{non}%
\end{equation}

\end{prop}

For any $q\in\left[  2,6\right]  $ we have a continuous embedding
$X\hookrightarrow L^{q}( \mathsf{R}^{3}) $. Consequently there is a constant
$\kappa_{q}>0$ such that%
\begin{equation}
\left\vert u\right\vert _{q}\leq\kappa_{q}\left\Vert u\right\Vert
\text{,\qquad for all }u\in X\text{.} \label{ka}%
\end{equation}

If $0<\lambda_{1}$, it is easy to see that $\Phi$ has the mountain pass
geometry. This case is simple and will be omitted here. For the proof of
Theorem \ref{t1}, since $0$ is not an eigenvalue of \eqref{eig}, we may assume
that $0\in(\lambda_{\ell},\lambda_{\ell+1})$ for some $\ell\geq1$. Let%
\begin{equation}
X^{-}=\operatorname*{span}\left\{  \phi_{1},\cdots,\phi_{\ell}\right\}
\text{,\qquad}X^{+}=(X^{-})^{\bot}\text{.} \label{xp}%
\end{equation}
Then $X^{-}$ and $X^{+}$ are the negative space and positive space of the
quadratic form%
\[
Q(u)=\frac{1}{2}\int_{\mathsf{R}^{3}}\left(  \left\vert \nabla u\right\vert
^{2}+V(x)u^{2}\right)  \mathrm{d}x
\]
respectively, note that $\dim X^{-}=\ell<\infty$. Moreover, there is a
positive constants $\kappa$ such that%
\begin{equation}
\pm Q(u)\geq\kappa\left\Vert u\right\Vert ^{2}\text{,\qquad}u\in X^{\pm
}\text{.} \label{ek}%
\end{equation}

As we have mentioned, because of the term involving $\phi_{u}$ in \eqref{e3},
our functional $\Phi$ does not satisfy the geometric assumption of the linking
theorem. In fact, choose $\phi\in X^{+}$ with $\left\Vert \phi\right\Vert =1$.
For $R>r>0$ set%
\[
N=\left\{  \left.  u\in X^{+}\right\vert \,\left\Vert u\right\Vert =r\right\}
\text{,\qquad}M=\left\{  \left.  u\in X^{-}\oplus\mathsf{R}^{+}\phi\right\vert
\,\left\Vert u\right\Vert \leq R\right\}  \text{,}%
\]
then $M$ is a submanifold of $X^{-}\oplus\mathsf{R}^{+}\phi$ with boundary
$\partial M$. We do have%
\[
b=\inf_{N}\Phi>0\text{,\qquad}\sup_{u\in\partial M,\left\Vert u\right\Vert
=R}\Phi<0
\]
provided $R$ is large enough. However, for $u\in X^{-}$ we have%
\[
\Phi(u)=Q(u)+\frac{1}{4}\int_{\mathsf{R}^{3}}\phi_{u}u^{2}\mathrm{d}%
x-\int_{\mathsf{R}^{3}}G(u)\mathsf{d}x\text{.}%
\]
Because $\phi_{u}\geq0$, the term involving $\phi_{u}$ may be very large and
for some point $u\in\partial M\cap X^{-}$ we may have $\Phi(u)>b$. Therefore
the following geometric assumption of the linking theorem
\[
b=\inf_{N}\Phi>\sup_{\partial M}\Phi
\]
can not be satisfied. Fortunately, we can apply the local linking theorem of
Luan and Mao \cite{MR2135818} (see also Li and Willem \cite{MR1312028}) to
overcome this difficulty and find critical points of $\Phi$.

Recall that by definition, $\Phi$ has a local linking at $0$ with respect to
the direct sum decomposition $X=X^{-}\oplus X^{+}$, if there is $\rho>0$ such
that%
\begin{equation}
\left\{
\begin{array}
[c]{ll}%
\Phi( u) \leq0\text{,} & \text{for }u\in X^{-}\text{, }\left\Vert u\right\Vert
\leq\rho\text{,}\\
\Phi( u) \ge0\text{,} & \text{for }u\in X^{+}\text{, }\left\Vert u\right\Vert
\leq\rho\text{.}%
\end{array}
\right.  \label{loc}%
\end{equation}
It is then clear that $0$ is a (trivial) critical point of $\Phi$. Next, we
consider two sequences of finite dimensional subspaces%
\[
X_{0}^{\pm}\subset X_{1}^{\pm}\subset\cdots\subset X^{\pm}%
\]
such that
\[
X^{\pm}=\overline{\bigcup_{n\in\mathsf{N}}X_{n}^{\pm}}\text{.}%
\]
For multi-index $\alpha=( \alpha^{-},\alpha^{+}) \in\mathsf{N}^{2}$ we set
$X_{\alpha}=X_{\alpha^{-}}^{-}\oplus X_{\alpha^{+}}^{+}$ and denote by
$\Phi_{\alpha}$ the restriction of $\Phi$ on $X_{\alpha}$. A sequence
$\left\{  \alpha_{n}\right\}  \subset\mathsf{N}^{2}$ is admissible if, for any
$\alpha\in\mathsf{N}^{2}$, there is $m\in\mathsf{N}$ such that $\alpha
\leq\alpha_{n}$ for $n\geq m$; where for $\alpha,\beta\in\mathsf{N}^{2}$,
$\alpha\leq\beta$ means $\alpha^{\pm}\leq\beta^{\pm}$. Obviously, if $\left\{
\alpha_{n}\right\}  $ is admissible, then any subsequence of $\left\{
\alpha_{n}\right\}  $ is also admissible.

\begin{defn}
[{\cite[Definition 2.2]{MR2135818}}]We say that $\Phi\in C^{1}(X)$ satisfies
the Cerami type condition $(C)^{\ast}$, if whenever $\left\{  \alpha
_{n}\right\}  \subset\mathsf{N}^{2}$ is admissible, any sequence $\left\{
u_{n}\right\}  \subset X$ such that%
\begin{equation}
u_{n}\in X_{\alpha_{n}}\text{,\qquad}\sup_{n}\Phi(u_{n})<\infty\text{,\qquad
}\left(  1+\left\Vert u_{n}\right\Vert \right)  \Vert\Phi_{\alpha_{n}}%
^{\prime}(u_{n})\Vert_{X_{\alpha_{n}}^{\ast}}\rightarrow0 \label{ps}%
\end{equation}
contains a subsequence which converges to a critical point of $\Phi$.
\end{defn}

\begin{thm}
[{Local Linking Theorem, \cite[Theorem 2.2]{MR2135818}}]\label{t4}Suppose that
$\Phi\in C^{1}(X)$ has a local lingking at $0$, $\Phi$ satisfies $(C)^{\ast}$,
$\Phi$ maps bounded sets into bounded sets and, for every $m\in\mathsf{N}$,%
\begin{equation}
\Phi(u)\rightarrow-\infty\text{,\qquad as }\left\Vert u\right\Vert
\rightarrow\infty\text{, }u\in X^{-}\oplus X_{m}^{+}\text{.} \label{ant}%
\end{equation}
Then $\Phi$ has a nontrivial critical point.
\end{thm}

\begin{rem}
Theorem \ref{t4} is a generalization of the well known local linking theorem
of Li and Willem \cite[Theorem 2]{MR1312028}, where instead of $\left(
C\right)  ^{\ast}$, the stronger Palais-Smale type condition $\left(
PS\right)  ^{\ast}$ is assumed.
\end{rem}

For the proof of Theorem \ref{t2} we will use the fountain theorem of Bartsch
\cite{MR1219237}, see also \cite[Theorem 3.6]{MR1400007}. For $k=1,2,\cdots$,
let%
\begin{equation}
Y_{k}=\operatorname*{span}\left\{  \phi_{1},\cdots,\phi_{k}\right\}
\text{,\qquad}Z_{k}=\overline{\operatorname*{span}\left\{  \phi_{k},\phi
_{k+1},\cdots\right\}  }\text{.} \label{dec}%
\end{equation}

\begin{thm}
[Fountain Theorem]\label{t3}Assume that the \emph{even} functional $\Phi\in
C^{1}( X) $ satisfies the $( PS) $ condition. If there exists $k_{0}>0$ such
that for $k\geq k_{0}$ there exist $\rho_{k}>r_{k}>0$ such that

\begin{enumerate}
\item[$( \mathrm{i}) $] $b_{k}=\inf\limits_{u\in Z_{k},\left\Vert u\right\Vert
=r_{k}}\Phi( u) \rightarrow+\infty$, as $k\rightarrow\infty$,

\item[$( \mathrm{ii}) $] $a_{k}=\max\limits_{u\in Y_{k},\left\Vert
u\right\Vert =\rho_{k}}\Phi( u) \leq0$,
\end{enumerate}
then $\Phi$ has a sequence of critical points $\left\{  u_{k}\right\}  $ such
that $\Phi( u_{k}) \rightarrow+\infty$.
\end{thm}

\section{Proofs of Theorems \ref{t1} and \ref{t2}}

To study the functional $\Phi$, it will be convenient to write it in a form in
which the quadratic part is $\left\Vert u\right\Vert ^{2}$. Let $f(t)=g(t)+mt$%
. Then, by a simple computation, we have%
\begin{equation}
F(t):=\int_{0}^{t}f\left(  \tau\right)  \mathrm{d}\tau\leq\frac{t}%
{4}f(t)+\frac{\tilde{b}}{4}t^{2}\text{,\qquad where }\tilde{b}=b+m>0\text{.}
\label{e4}%
\end{equation}
Note that by \eqref{e2} we easily have%
\begin{equation}
\lim_{\left\vert t\right\vert \rightarrow\infty}\frac{f(t)}{t^{3}}%
=+\infty\text{.} \label{e7}%
\end{equation}
Moreover, using \eqref{e} we get%
\[
\lim_{\left\vert t\right\vert \rightarrow0}\frac{f(t)t}{t^{4}}=\lim
_{\left\vert t\right\vert \rightarrow0}\left(  \frac{t^{2}}{t^{4}}\cdot
\frac{g(t)t+mt^{2}}{t^{2}}\right)  =+\infty\text{.}%
\]
Therefore, there is $\Lambda>0$ such that%
\begin{equation}
f(t)t\geq-\Lambda t^{4}\text{,\qquad all }t\in\mathsf{R}\text{.} \label{e71}%
\end{equation}

With the modified nonlinearity $f$, our functional $\Phi:X\rightarrow
\mathsf{R}$ can be rewritten as follows:%
\begin{equation}
\Phi(u)=\frac{1}{2}\left\Vert u\right\Vert ^{2}+\frac{1}{4}\int_{\mathsf{R}%
^{3}}\phi_{u}u^{2}\mathrm{d}x-\int_{\mathsf{R}^{3}}F(u)\mathrm{d}x\text{.}
\label{ab}%
\end{equation}
Note that this does not imply that $\Phi$ has the mountain pass geometry,
because unlike in \eqref{e3}, as $\left\Vert u\right\Vert \rightarrow0$ the
last term in \eqref{ab} is not $o(\left\Vert u\right\Vert ^{2})$ anymore. The
derivative of $\Phi$ is given below:%
\[
\langle\Phi^{\prime}(u),v\rangle=\left\langle u,v\right\rangle +\int
_{\mathsf{R}^{3}}\phi_{u}uv\mathrm{d}x-\int_{\mathsf{R}^{3}}f(u)v\mathrm{d}%
x\text{.}%
\]

\begin{lem}
\label{lemma}Suppose $(V)$, $(g_{0})$ and $(g_{1})$ are satisfied, then $\Phi$
satisfies the $(C)^{\ast}$ condition.
\end{lem}

\begin{pf}
Suppose $\left\{  u_{n}\right\}  $ is a sequence satisfying \eqref{ps}, where
$\left\{  \alpha_{n}\right\}  \subset\mathsf{N}^{2}$ is admissible. We must
prove that $\left\{  u_{n}\right\}  $ is bounded.

We may assume $\left\Vert u_{n}\right\Vert \rightarrow\infty$ for a
contradiction. By \eqref{ps} and noting that%
\begin{equation}
\langle\Phi_{\alpha_{n}}^{\prime}(u_{n}),u_{n}\rangle=\langle\Phi^{\prime
}(u_{n}),u_{n}\rangle\label{res}%
\end{equation}
since $u_{n}\in X_{\alpha_{n}}$, for large $n$, using \eqref{e4} we have%
\begin{align}
4\cdot\sup_{n}\Phi(u_{n})+\left\Vert u_{n}\right\Vert  &  \geq4\Phi
(u_{n})-\langle\Phi_{\alpha_{n}}^{\prime}(u_{n}),u_{n}\rangle\nonumber\\
&  =\left\Vert u_{n}\right\Vert ^{2}+\int_{\mathsf{R}^{3}}\left(
f(u_{n})u_{n}-4F(u_{n})\right)  \mathrm{d}x\nonumber\\
&  \geq\left\Vert u_{n}\right\Vert ^{2}-\tilde{b}\int_{\mathsf{R}^{3}}%
u_{n}^{2}\mathrm{d}x\text{.} \label{e6}%
\end{align}
Let $v_{n}=\left\Vert u_{n}\right\Vert ^{-1}u_{n}$. Up to a subsequence, by
the compact embedding $X\hookrightarrow L^{2}(\mathsf{R}^{3})$ we deduce%
\[
v_{n}\rightharpoonup v\text{ in }X\text{,\qquad}v_{n}\rightarrow v\text{ in
}L^{2}(\mathsf{R}^{3})\text{.}%
\]
Multiplying by $\left\Vert u_{n}\right\Vert ^{-2}$ to both sides of \eqref{e6}
and letting $n\rightarrow\infty$, we obtain%
\[
\tilde{b}\int_{\mathsf{R}^{3}}v^{2}\mathrm{d}x\geq1\text{.}%
\]
Consequently, $v\neq0$.

Using \eqref{e71} and \eqref{ka} with $q=4$, we have%
\begin{align}
\int_{v=0}\frac{f(u_{n})u_{n}}{\left\Vert u_{n}\right\Vert ^{4}}\mathrm{d}x
&  =\int_{v=0}\frac{f(u_{n})u_{n}}{u_{n}^{4}}v_{n}^{4}\mathrm{d}x\nonumber\\
&  \geq-\Lambda\int_{v=0}v_{n}^{4}\mathrm{d}x\nonumber\\
&  \geq-\Lambda\int_{\mathsf{R}^{3}}v_{n}^{4}\mathrm{d}x=-\Lambda\left\vert
v_{n}\right\vert _{4}^{4}\geq-\Lambda\kappa_{4}^{4}>-\infty\text{.} \label{ev}%
\end{align}
For $x\in\left\{  x\in\mathsf{R}^{3}\mid v\neq0\right\}  $, we have
$\left\vert u_{n}(x)\right\vert \rightarrow+\infty$. By \eqref{e7} we get%
\begin{equation}
\frac{f(u_{n}(x))u_{n}(x)}{\left\Vert u_{n}\right\Vert ^{4}}=\frac
{f(u_{n}(x))u_{n}(x)}{u_{n}^{4}(x)}v_{n}^{4}(x)\rightarrow+\infty\text{.}
\label{inff}%
\end{equation}
Consequently, using \eqref{ev}, \eqref{inff} and the Fatou lemma we obtain%
\begin{equation}
\int_{\mathsf{R}^{3}}\frac{f(u_{n})u_{n}}{\left\Vert u_{n}\right\Vert ^{4}%
}\mathrm{d}x\geq\int_{v\neq0}\frac{f(u_{n})u_{n}}{\left\Vert u_{n}\right\Vert
^{4}}\mathrm{d}x-\Lambda\kappa_{4}^{4}\rightarrow+\infty\text{.} \label{xm}%
\end{equation}
Since $\left\{  u_{n}\right\}  $ is a sequence satisfying \eqref{ps}, using
\eqref{non} and \eqref{xm}, for large $n$ we obtain
\begin{align}
4a_{1}+1  &  \geq\frac{1}{\left\Vert u_{n}\right\Vert ^{4}}\left(  \left\Vert
u_{n}\right\Vert ^{2}+\int_{\mathsf{R}^{3}}\phi_{u_{n}}u_{n}^{2}%
\mathrm{d}x-\left\langle \Phi^{\prime}(u_{n}),u_{n}\right\rangle \right)
\nonumber\\
&  =\int_{\mathsf{R}^{3}}\frac{f(u_{n})u_{n}}{\left\Vert u_{n}\right\Vert
^{4}}\mathrm{d}x\rightarrow+\infty\text{,} \label{inf}%
\end{align}
a contradiction.

Therefore, $\left\{  u_{n}\right\}  $ is bounded in $X$. Now, by the argument
of \cite[pp. 4933]{MR2548724} and using \eqref{res}, the compact embedding
$X\hookrightarrow L^{2}(\mathsf{R}^{3})$ and%
\[
X=\overline{\bigcup_{n\in\mathsf{N}}X_{\alpha_{n}}}\text{,}%
\]
we can easily prove that $\left\{  u_{n}\right\}  $ has a subsequence
converging to a critical point of $\Phi$.
\end{pf}

\begin{lem}
\label{lem1}Suppose $( V) $, $( g_{0}) $ and $( g_{1}) $ are satisfied, then
$\Phi$ satisfies the $( PS) $ condition.
\end{lem}

\begin{pf}
Under the assumption there exists $\tilde{\Lambda}>0$ such that%
\begin{equation}
F(t)\geq-\tilde{\Lambda}t^{4}\text{,\qquad}\lim_{\left\vert t\right\vert
\rightarrow\infty}\frac{F(t)}{t^{4}}=+\infty\text{.} \label{ff}%
\end{equation}
Let $\left\{  u_{n}\right\}  $ be a $\left(  PS\right)  $ sequence, that is
$\sup_{n}\left\vert \Phi(u_{n})\right\vert <\infty$, $\Phi^{\prime}%
(u_{n})\rightarrow0$. We only need to show that $\left\{  u_{n}\right\}  $ is
bounded. If $\left\{  u_{n}\right\}  $ is not bounded, similar to the first
part in the proof of Lemma \ref{lemma}, we may assume that for some $v\neq0$,%
\[
\left\Vert u_{n}\right\Vert ^{-1}u_{n}\rightharpoonup v\text{\qquad in
}X\text{.}%
\]
Since $v\neq0$, similar to \eqref{xm}, using \eqref{ff} we deduce%
\[
\int_{\mathsf{R}^{3}}\frac{F(u_{n})}{\left\Vert u_{n}\right\Vert ^{4}%
}\mathrm{d}x\rightarrow+\infty\text{.}%
\]
Therefore%
\begin{align*}
a_{1}+1  &  \geq\frac{1}{\left\Vert u_{n}\right\Vert ^{4}}\left(  \frac{1}%
{2}\left\Vert u_{n}\right\Vert ^{2}+\frac{1}{4}\int_{\mathsf{R}^{3}}%
\phi_{u_{n}}u_{n}^{2}\mathrm{d}x-\Phi(u_{n})\right) \\
&  =\int_{\mathsf{R}^{3}}\frac{F(u_{n})}{\left\Vert u_{n}\right\Vert ^{4}%
}\mathrm{d}x\rightarrow+\infty\text{,}%
\end{align*}
a controdiction.
\end{pf}

\begin{lem}
\label{lem5}Under the assumptions $( V) $, $( g_{0}) $, the functional $\Phi$
has a local linking at $0$ with respect to the decomposition $X=X^{-}\oplus
X^{+}$.
\end{lem}

\begin{pf}
By $(g_{0})$, there exists $C>0$ such that%
\begin{equation}
\left\vert G(u)\right\vert \leq\frac{\kappa}{2\kappa_{2}^{2}}\left\vert
u\right\vert ^{2}+C\kappa\left\vert u\right\vert ^{p}\text{.} \label{inq}%
\end{equation}
If $u\in X^{-}$, then using \eqref{non} and \eqref{ek} we deduce%
\begin{align}
\Phi(u)  &  =Q(u)+\frac14\int_{\mathsf{R}^{3}}\phi_{u}u^{2}\mathrm{d}%
x-\int_{\mathsf{R}^{3}}G(u)\mathrm{d}x\nonumber\\
&  \leq-\kappa\left\Vert u\right\Vert ^{2}+a_{1}\left\Vert u\right\Vert
^{4}+\frac{\kappa}{2\kappa_{2}^{2}}\left\vert u\right\vert _{2}^{2}%
+C\kappa\left\vert u\right\vert _{p}^{p}\nonumber\\
&  \leq-\frac{\kappa}{2}\left\Vert u\right\Vert ^{2}+a_{1}\left\Vert
u\right\Vert ^{4}+C_{1}\left\Vert u\right\Vert ^{p}\text{,} \label{lo1}%
\end{align}
where $C_{1}=C\kappa\kappa_{p}^{p}$. Similarly, for $u\in X^{+}$ we have%
\begin{equation}
\Phi(u)\geq\frac{\kappa}{2}\left\Vert u\right\Vert ^{2}-C_{1}\left\Vert
u\right\Vert ^{p}\text{.} \label{lo2}%
\end{equation}
Since $p>4$, the desired result \eqref{loc} follows from \eqref{lo1} and \eqref{lo2}.
\end{pf}

\begin{lem}
\label{lem6}Let $Y$ be a finite dimensional subspace of $X$, then $\Phi$ is
anti-coercive on $Y$, that is%
\[
\Phi( u) \rightarrow-\infty\text{,\qquad as }\left\Vert u\right\Vert
\rightarrow\infty\text{, }u\in Y\text{.}%
\]

\end{lem}

\begin{pf}
If the conclusion is not true, we can choose $\left\{  u_{n}\right\}  \subset
Y$ and $\beta\in\mathsf{R}$ such that%
\begin{equation}
\left\Vert u_{n}\right\Vert \rightarrow\infty\text{,\qquad}\Phi(u_{n}%
)\geq\beta\text{.} \label{con}%
\end{equation}
Let $v_{n}=\left\Vert u_{n}\right\Vert ^{-1}u_{n}$. Since $\dim Y<\infty$, up
to a subsequence we have%
\[
\left\Vert v_{n}-v\right\Vert \rightarrow0\text{,\qquad}v_{n}(x)\rightarrow
v(x)\text{ a.e. }\mathsf{R}^{3}%
\]
for some $v\in Y$, with $\Vert v\Vert=1$. Since $v\neq0$, similar to
\eqref{xm}, using \eqref{ff} we have%
\begin{equation}
\int_{\mathsf{R}^{3}}\frac{F(u_{n})}{\left\Vert u_{n}\right\Vert ^{4}%
}\mathrm{d}x\rightarrow+\infty\text{.} \label{xx}%
\end{equation}
Using \eqref{non} and \eqref{xx} we deduce
\[
\Phi(u_{n})=\left\Vert u_{n}\right\Vert ^{4}\left(  \frac{1}{2\left\Vert
u_{n}\right\Vert ^{2}}+\frac{1}{4\left\Vert u_{n}\right\Vert ^{4}}%
\int_{\mathsf{R}^{3}}\phi_{n}u_{n}^{2}\mathrm{d}x-\int_{\mathsf{R}^{3}}%
\frac{F(u_{n})}{\left\Vert u_{n}\right\Vert ^{4}}\mathrm{d}x\right)
\rightarrow-\infty\text{,}%
\]
a contradiction with \eqref{con}.
\end{pf}

Now, we are ready to prove our main results.

\begin{pf}
[Proof of Theorem \ref{t1}]We will find a nontrivial critical point of $\Phi$
via Theorem \ref{t4}. Using Proposition \ref{p1}, it is easy to see that
$\Phi$ maps bounded sets into bounded sets. In Lemmas \ref{lemma} and
\ref{lem5} we see that $\Phi$ satisfies $(C)^{\ast}$ and $\Phi$ has a local
linking at $0$. It suffices to verify \eqref{ant}. Since $\dim(X^{-}\oplus
X_{m}^{+})<\infty$, this is a consequence of Lemma \ref{lem6}.
\end{pf}

\begin{pf}
[Proof of Theorem \ref{t2}]By Remark \ref{re13} (i), it suffices to find a
sequence of critical points $\{u_{n}\}$ of $\Phi$ such that $\Phi(u_{n}%
)\to+\infty$.

We define subspaces $Y_{k}$ and $Z_{k}$ of $X$ as in \eqref{dec}. Since $g$ is
odd, $\Phi$ is an even functional. By Lemma \ref{lem1} we know that $\Phi$
satisfies $( PS) $. It suffices to verify (i) and (ii) of Theorem \ref{t3}.

\emph{Verification of (i)}. We assume that $0\in\lbrack\lambda_{\ell}%
,\lambda_{\ell+1})$. Then if $k>\ell$ we have $Z_{k}\subset X^{+}$, where
$X^{+}$ is defined in \eqref{xp}. Now, by \eqref{ek}, we have%
\begin{equation}
Q(u)\geq\kappa\left\Vert u\right\Vert ^{2}\text{,\qquad}u\in Z_{k}\text{.}
\label{q}%
\end{equation}
We recall that by \cite[Lemma 2.5]{MR2548724},%
\begin{equation}
\beta_{k}=\sup_{u\in Z_{k},\left\Vert u\right\Vert =1}\left\vert u\right\vert
_{p}\rightarrow0\text{,\qquad as }k\rightarrow\infty\text{.} \label{bt}%
\end{equation}
Let $r_{k}=(Cp\beta_{k}^{p})^{1/(2-p)}$, where $C$ is chosen in \eqref{inq}.
For $u\in Z_{k}$ with $\left\Vert u\right\Vert =r_{k}$, using \eqref{inq} and
\eqref{q} we deduce%
\begin{align*}
\Phi(u)  &  =Q(u)+\frac{1}{4}\int_{\mathsf{R}^{3}}\phi_{u}u^{2}\mathrm{d}%
x-\int_{\mathsf{R}^{3}}G(u)\mathrm{d}x\\
&  \geq\kappa\left\Vert u\right\Vert ^{2}-\frac{\kappa}{2\kappa_{2}^{2}%
}\left\vert u\right\vert _{2}^{2}-C\kappa\left\vert u\right\vert _{p}^{p}\\
&  \geq\kappa\left(  \frac{1}{2}\left\Vert u\right\Vert ^{2}-C\beta_{k}%
^{p}\left\Vert u\right\Vert ^{p}\right)  =\kappa\left(  \frac{1}{2}-\frac
{1}{p}\right)  \left(  Cp\beta_{k}^{p}\right)  ^{2/(2-p)}\text{.}%
\end{align*}
Since $\beta_{k}\rightarrow0$ and $p>2$, it follows that
\[
b_{k}=\inf\limits_{u\in Z_{k},\left\Vert u\right\Vert =r_{k}}\Phi( u)
\rightarrow+\infty\text{.}
\]

\emph{Verification of (ii)}. Since $\dim Y_{k}<\infty$, this is a consequence
of Lemma \ref{lem6}.
\end{pf}

\begin{rem}
From the proof of Lemmas \ref{lem1} and \ref{lem6}, it is easy to see that if
in $(g_{1})$ we replace \eqref{e2} by the following weaker condition%
\[
\lim_{\left\vert t\right\vert \rightarrow\infty}\frac{F(t)}{t^{4}}%
=+\infty\text{,}%
\]
then Theorem \ref{t2} remains valid.
\end{rem}

\end{document}